\newcommand{\Comments}[1]{\\ \fbox{\fbox{\rm #1}}\\}
\renewcommand{\Comments}[1]{\mbox{}}
\documentclass[11pt]{article}
 \usepackage{amsmath,amsfonts,amssymb}
 \usepackage{amsfonts}
 \usepackage{amssymb}
\newcommand{\ls}[1]
   {\dimen0=\fontdimen6\the\font \lineskip=#1\dimen0
 \advance\lineskip.5\fontdimen5\the\font \advance\lineskip-\dimen0
 \lineskiplimit=.9\lineskip \baselineskip=\lineskip
 \advance\baselineskip\dimen0 \normallineskip\lineskip
 \normallineskiplimit\lineskiplimit \normalbaselineskip\baselineskip
 \ignorespaces }
 \numberwithin{equation}{section}
 \newcommand{\beq}{\begin{equation}}
 \newcommand{\eeq}{\end{equation}}
 \newcommand{\beqa}{\begin{eqnarray}}
 \newcommand{\eeqa}{\end{eqnarray}}
 \newcommand{\beqann}{\begin{eqnarray*}}
 \newcommand{\eeqann}{\end{eqnarray*}}

 \newcommand{\pf}{\noindent \mbox{{\bf Proof}: }}
 \def\squarebox#1{\hbox to #1{\hfill\vbox to #1{\vfill}}}
 
 \newcommand{\grad}[1]{\nabla #1}

 \newcommand{\hilbert}{\bigcirc\kern -0.8em{\rm\scriptstyle {H}\;}}
 \newcommand{\FFF}{F}

 \newcommand{\calS}{{\mathcal S}}%
 \newcommand{\pair}[4]{\raisebox{-1.2ex}{\mbox{\tiny$#1$}}\!
       \langle#2,#3\rangle\!\raisebox{-1.2ex}{\mbox{\tiny $#4$}}}
 
 \newcommand{\clfs}[2]{\pair{W^*}{#2}{#1}{\Wss}}

 \newcommand{\ddom}{{\rm \bf{dom}}}

 \newcommand{\ddomdel}[1]{\ddom_{#1}\ddelta}

 \newcommand{\ddelta}{\delta\hspace{-.2cm}\delta}

 \newcommand{\Wss}{W^{**}}
 
 \newcommand{\KKK}{\mathbf K}

\newcommand{\tr}{{\rm tr}}
\begin{document}
 \noindent
 {\em Erratum}\\

 \noindent
 {\bf \large Correction to ``The divergence of Banach space valued
 random variables on Wiener space", Prob. Th. Rel. Fields 132, 291-320 (2005)}\\

 E. Mayer-Wolf \footnote{Department of Mathematics, Technion I.I.T., Haifa, Israel}
       \ \ and M. Zakai\footnote{Department of Electrical Engineering, Technion I.I.T., Haifa, Israel}

 \vspace{1cm}

 \noindent
 We are grateful to J. Maas and J. Van Neerven  for drawing our attention to the two mistakes addressed below.\\

 \noindent
 Corollary~3.5 (as well as Corollary~3.17a whose proof relies on it) should be ignored since
 the inequality $|F_n|_{p,1}\!\le\!\|F_n\|_{p,1}$ in its proof is false, and we have been unable
 to find a simple alternative argument.\\

 \noindent
 More importantly, in Proposition~3.14 and Proposition~3.18 one needs to add the assumption

 \noindent
 {\bf (A)}\ \ \ $Y^{**}$ has the Radon Nykodim property (RNP) with respect to $\mu$\ \ (cf. \cite{DU})

 \noindent
 on the Banach space $Y$, as we shall now explain. Note (Section III.3 in \cite{DU}) that $Y^{**}$ has the RNP
 with respect to any measure if, for example, $Y^{**}$ is separable or $Y$ is reflexive.

 \noindent
 Unfortunately, in the proof of Proposition~3.14, the natural imbedding of\
 $L^p(\mu;Y^{**})$\ in the operator space \ $L(Y^*,L^p(\mu))$\ was erroneously
 claimed to be surjective. In addition, the observation associated
 with (3.12) was also incorrect as stated (although not used in the rest of
 the paper). We restate this observation in {\bf (1)} below, and prove it under the
 additional assumption {\bf (A)}; it will then be used in {\bf (2)} to replace the incorrect 
 proof of Proposition~3.14.\\
 \begin{itemize}
  \item[{\bf (1)}] Assume {\bf (A)} and let $\KKK\!\in\!L^p(\mu;L(W^*,Y^{**})$.\ \
   If~(3.12) holds for some $\gamma\!>\!0$ and all $F\!\in{\mathcal S}(Y^*)$\
   then $\KKK\!\in\!\ddomdel{p,Y^{**}}$\ \ \ (the converse is obvious).

   \pf The bound~(3.12) implies the existence of a $\Lambda_{\KKK}\!\in\!L^q(\mu,Y^*)^*$ such that
   \[ E\tr \left(\KKK^T\grad^{W^*}F\right)=\Lambda_{\KKK}(F)
                    \hspace{.75 cm} \forall F\!\in\!{\mathcal S}(Y^*)\ .\]
   Due to assumption~{\bf (A)},\ $L^q(\mu,Y^*)^*$ can be
   identified with\ $L^p(\mu,Y^{**})$\ (cf. Theorem IV.1 in \cite{DU})\ in the sense that there
   exists a $\ddelta\KKK\!\in\!L^p(\mu,Y^{**})$\ such
   that\ $\Lambda_{\KKK}(F)$ is given by \ \
   $E\pair{Y^*}{F}{\ddelta\KKK}{Y^{**}}$. Thus, for any  $F\!\in\!\calS(Y^{*})$
   \[ \hspace*{3cm}E\tr \left(\KKK^T\grad^{W^*}F\right)=E\pair{Y^{**}}{\ddelta\KKK}{F}{Y^{***}}
    \ .\hspace{3cm}(\dag)\]
   For\ \ $F=\phi(\delta(e_1),\ldots,\delta(e_m))\otimes l$,\ \ \ \
   ($\{e_1,\ldots,e_m\}\!\subset\!W^*$, orthonormal in $H$,\ and $l\!\in\!Y^*$),\ \ \
   $(\dag)$\ amounts to
    $\pair{Y^{**}}{E\left(\sum_{i=1}^m\partial_i\phi\,Ke_i\right)}{l}{Y^{***}}
      \!\!=\!\!\pair{Y^{**}}{E\phi\,\ddelta\KKK}{l}{Y^{***}}$
   which, if true $\forall l\!\in\!Y^*$, is true $\forall l\!\in\!Y^{***}$ as well.
   Thus\ $(\dag)$ holds for all $F\!\in\!\calS(Y^{***})$, which means that
   $\KKK\!\in\!\ddomdel{p,Y^{**}}$.
  \item[{\bf (2)}] We now present a modified proof of the ``if" implication in
   the first statement of Proposition 3.14, using the characterization provided
   by {\bf (1)} instead of the erroneous identification of\ $L^p(\mu;Y^{**})$\
   and\ $L(Y^*,L^p(\mu))$\ mentioned above:\\

   \noindent
   It follows from~(3.13) that there exists a $\Delta_{_\KKK}\!\in\!L(Y^*,L^p(\mu))$
   such that for all $l\!\in\!Y^*$
 \[    \hspace*{5.5cm}\ddelta(\KKK^Tl)=\Delta_{_\KKK}(l)\hspace{5.5cm} (3.17)\]
   so that, for any\ $\FFF=\sum_{j=1}^m \Phi_jl_j\ \ \in\!\calS(Y^*)$\
   \begin{eqnarray*}
    E{\rm tr}\,(\KKK^T\grad{\FFF})
      &=&E\sum_{j=1}^m {\rm tr}\,\KKK^T\grad(\Phi_j\,l_j)
         =\sum_{j=1}^m E\clfs{\KKK^Tl_j}{\grad{\Phi_j}}\\
      &=&\sum_{j=1}^m E\delta(\KKK^Tl_j)\Phi_j
         =\sum_{j=1}^m E\Delta_{_\KKK}(l_j)\Phi_j\\
      &=&E\Delta_{_\KKK}\!\left(\!\sum_{j=1}^m\Phi_jl_j\!\right)=E\Delta_{_\KKK}(F)
  \end{eqnarray*}
  and thus for any $q\!\ge\!1$,\ and with $\|\!\Delta_{_\KKK}\!\|$\ denoting the operator norm,
   \[ \left|\,E{\rm tr}\,(\KKK^T\grad{\FFF})\,\right|\! \le\!E\|\!\Delta_{_\KKK}\!\|\|F\|_{_{Y^*}}\!
   \le\!\|\!\Delta_{_\KKK}\!\|\left(E\|F\|^q_{_{Y^*}}\right)^{1/q}\]
   so that from {\bf (1)}\ it follows that $\KKK\!\in\!\ddomdel{p,Y^{**}}$.
 \end{itemize}
 

\begin{thebibliography}{99}
 \bibitem{DU} J. Diestel and J.J. Uhl, Jr, {\em Vector Measures},\
  AMS Math. Surv. 15 (1977)
 \end{thebibliography}
 \end{document}